\documentclass[11pt]{article}

\usepackage{amsmath}
\usepackage{amssymb}
\usepackage{amsthm}
\usepackage{color}
\usepackage{latexsym}
\usepackage[english]{babel}

\newtheorem{Thm}{Theorem}

\newtheorem{Prop}[Thm]{Proposition}
\newtheorem{Lem}[Thm]{Lemma}
\newtheorem{Def}[Thm]{Definition}
\newtheorem{rem}{Remark}
\newcommand{\dn}[1]{\dfrac{\partial #1}{\partial \nu}}

\newcommand{\dns}[2]{\dfrac{\partial #1}{\partial \nu_{#2}}}
\newcommand\E{\varepsilon}

\newcommand{\Bi}{B_{\varepsilon}(x_1)}

\newcommand{\fun}[3]{#1 : #2 \longrightarrow #3}
\newcommand{\dem}{\textbf{Proof :}}
\newcommand {\QED}{\hfill{\rule{2mm}{2mm}}}
\newcommand{\Gm}{\Gamma}
\newcommand{\tn}{\textnormal}
\newcommand{\ds}{\displaystyle}

\def\Sg{\Sigma}
\def\la{\lambda}
\def\Om{\Omega}

\def\RR{\mathbb{R}}
\def\CC{\mathbb{C}}

\def\1{1 \!\! l}

\def\Div{\textnormal{div}}

\title{An inverse problem  for Schr\"odinger equations with discontinuous main  
coefficient.} 

\date{}

\author{Lucie Baudouin\footnote{e-mail: {\tt baudouin@laas.fr}}\\
{\it\footnotesize LAAS - CNRS; Universit\'e de Toulouse; 7, avenue du Colonel Roche, F-31077 Toulouse, France.}\\ 
Alberto Mercado \footnote{e-mail: {\tt amercado@dim.uchile.cl}}\\
{\it\footnotesize  DIM-CMM, Universidad de Chile, Casilla 170/3 - Correo 3, Santiago, Chile.}}

\begin{document}

\maketitle

\begin{abstract}
This paper concerns  the inverse problem of retrieving a stationary potential for the Schr\"odinger evolution equation in a bounded domain of $\RR^N$ with Dirichlet data and discontinuous principal coefficient $a(x)$ from a single time-dependent Neumann boundary measurement.  We consider that the discontinuity of $a$ is located on a simple closed hyper-surface   called the {\it interface}, and $a$ is constant in each one of the interior and exterior domains with respect to this interface.  We prove uniqueness and lipschitz stability for this inverse problem under certain convexity hypothesis on the geometry of the interior domain and on the sign of the jump of $a$ at the interface. The proof is based on a global Carleman 
inequality for the  Schr\"odinger  equation with discontinuous coefficients, result also interesting by itself. 
\end{abstract}



\section{Introduction}\label{A1}

The method of Carleman estimates was introduced in the field of inverse problems by \textsc{Bukhgeim} and \textsc{Klibanov} in reference \cite{B-K} (see also \cite{B} and \cite{Kl}). The first known results concern uniqueness of inverse problems. Then, one of the first stability result for a multidimensional inverse problem, dealing with an hyperbolic equation, can be read in \cite{P-Y3} and is based on a modification of an idea of \cite{B-K}.

Carleman estimates techniques are presented in  \cite{Kliba-Timo} for standard coefficient inverse problems  for both linear and nonlinear partial differential equations;  one can also read in this book the construction of globally convergent numerical methods for coefficient inverse problems  and some concrete applied problems in geophysics, medical imaging and computational time reversal.

It is possible to obtain local Lipschitz stability around the single known solution, provided that this solution is regular enough and contains enough information (see \cite{K-M} and \cite{Kl}).
Actually, many of the results using the same strategy we can refer to concern the wave equation.
A complete list is too long to be given here but to cite some of them, related to the same kind of inverse problems of determining a potential and also using local or global Carleman estimates, see \cite{P-Y1} and \cite{Y} for a Dirichlet boundary data and a Neumann measurement and \cite{I-Y2} for a Neumann boundary data and a Dirichlet measurement.

Recently, global Carleman estimates and applications to one-measurement inverse problems were obtained in the case of variable but still regular coefficients, see \cite{I-Y1} for the isotropic case, and \cite{L-T-Y} and \cite{Bella} for the anisotropic case. It is interesting to note that these authors require a bound on the gradient of the coefficients, so that the idea of approximating discontinuous coefficients by smooth ones is not useful. Nevertheless, uniqueness and Lipschitz stability 
are obtained in \cite{B-M-O} for the inverse problem of retrieving a stationary potential for the 
wave equation with Dirichlet data and discontinuous principal coefficient from a 
single time-dependent Neumann boundary measurement. 

One can also note that a global Carleman estimate has been obtained  \cite{D-O-P}  for the heat equation with discontinuous coefficients. That work was initially motivated by the study of the exact null controllability of the semilinear heat equation, but the estimate has been recently used (see  \cite{Bella2} )  to prove local Lipschitz stability for a one measurement inverse problem. In this field, one should also read the recent works \cite{B-D-LR},  \cite{B-G-LR} and \cite{Pois}. 
Up to our knowledge, the result of determination of a time independent 
potential in Schr\"odinger evolution equation with discontinuous principal coefficient from a single time dependent measurement on the boundary is new. Concerning the simpler case of a ``classical" Schr\"odinger equation (with $a   =  1 $), one can have a look at  \cite{B-P}, where the Carleman estimate and the proof of the stability of the same inverse problem are maybe easier to read and the philosophy is the same.
For the same equation, one can find in \cite{M-O-R} a method with weight functions satisfying a relaxed pseudoconvexity condition, which allows to prove  Carleman inequalities with less restrictive boundary observations than in \cite{B-P}.
The authors of \cite{La-Tr-Zh}  deal  with Carleman estimates  for the Sch\"odinger equation with variable  (but regular) principal coefficient  and applications to controllability.
Let us notice that in the different context of Cauchy problem,  V. \textsc{Isakov} in \cite{I} uses local Carleman 
estimates for the Schr\"odinger equation to prove uniqueness of the solution.
Finaly, for the Schr\"odinger operator $i\partial_t + ~\Div (c\nabla)$ in an unbounded strip in $\RR^2$, reference \cite{C-C-G} gives a stability result for the diffusion coefficient $c$ in $H^1$ with only one observation in an unbounded domain. One will see in the proof of our main tool (an appropriate Carleman estimate) that it is based on the same strong pseudoconvexity condition $(H_4)$ for the weight $\psi$.

\subsection{Statement of the problem and main results.}

Let $T>0$ and let $\Omega \subset \mathbb{R}^N$ ($N \ge 2$) be a
bounded domain with $C^2$-boundary $\partial\Omega$. 
Throughout this paper, we use the following notations :
\begin{eqnarray*}
&&~\nabla v = \left(\frac {\partial v}{\partial x_1},\dots,
\frac {\partial v}{\partial x_N}\right),~\Delta v = \sum_{i=1}^{N}\frac {\partial^2 v}{\partial x_i^2}~,\\
&&~v'=\frac {\partial v}{\partial t}~\tn{ and }~v''=\frac {\partial^2 v}{\partial t^2}~,\\
&&~\nu\in\mathbb{R}^N~\tn{ denotes the unit outward normal vector to }
\partial\Omega,\\
&&~\frac {\partial v}{\partial \nu}=\nabla v.\nu~\tn{  is the normal
derivative.}\\
\end{eqnarray*}

We will work with the following Schr\"odinger equation :
\begin{equation}\label{Anonlin}
\left\lbrace
\begin{array}{ll}
iy'(x,t)+ ~\Div (a(x) \nabla y(x,t))+p(x)y(x,t)=0,&x\in\Omega,~t\in (0,T)\\
y(x,t)=h(x,t),&x \in\partial\Omega,~t\in (0,T)\\
y(x,0)=y_0(x),&x\in\Omega.
\end{array} \right.\\
\end{equation}

We consider in this paper the inverse problem of the determination of the coefficient $p$ of the lower
order term in Schr\"odinger equation (\ref{Anonlin}) from a single time dependent observation of
Neumann data $\frac{\partial y}{\partial\nu}$ on the boundary.

The major novelty of this paper is that we deal with a Schr\"odinger equation in a bounded domain of $\RR^N$  with discontinuous principal coefficient.
Indeed, let  $\Om$ and  $\Om_1$ be two open subsets of $\RR^N$ with smooth boundaries $\Gm$ and $\Gm_1$. We choose $\Om_1$ simply connected and  such that $\overline{\Om}_1 \subset \Om$ and we set $\Om_2 = \Om \backslash \overline{\Om}_1$. Thus, we have $\partial \Om_2 = \Gm \cup \Gm_1$ and we also set:
\begin{equation} 
a(x) = \left\{ \begin{array}{cc}
a_1    &x\in \Om_1 \\
a_2    &x\in \Om_2
\end{array}
\right. \notag
\end{equation}
with $a_j > 0$ for $j = 1,2$.  \\

Considering equation (\ref{Anonlin}), we know that for each 
$p \in L^{\infty}(\Om)$, $y_0 \in L^2(\Om)$ and $h \in L^2(\Gm\times(0,T))$, 
there exists a unique weak solution $y$ such that
$$y \in C([0,T];H^{-1}(\Om))\cap H^{-1}(0,T;L^2(\Om)).$$
The proof is based on a transposition method, as one can read in \cite{M}. Let us also notice that the regularity of $y$ also implies $\dn{y} \in H^{-2}(0,T;H^{-\frac32}(\Gm))$. \\

We will prove the well-posedness of the inverse problem consisting in retrieving the potential $p$ involved in equation (\ref{Anonlin}), knowing the flux (the normal derivative) of the solution $y(p)$ of (\ref{Anonlin}) on the boundary. It means that we will prove uniqueness and stability of the nonlinear inverse problem characterized by the  nonlinear application 
\begin{equation} \label {nonlinear}
\left.p\right|_{\Omega}\longmapsto\left.a_2\frac{\partial y}{\partial \nu}\right|_{\Gamma\times(0,T)}.
\end{equation}
We will more precisely answer the following questions.

\noindent\textbf{Uniqueness} : \\
Does the equality $~\dfrac{\partial y(p)}{\partial\nu}=\dfrac{\partial y(q)}{\partial \nu}~$ on $~\Gamma\times(0,T)~$ 
imply $~p=q~$ on $~\Omega$ ?\\

\noindent\textbf{Stability} : \\
Is it possible to estimate $~\left.q-p\right|_{\Omega}~$ by  $
~\left.\dfrac{\partial y(q)}{\partial \nu}-\dfrac{\partial y(p)}{\partial \nu}\right|_{\Gamma\times(0,T)}~ $ in suitable norms~?\\

Indeed, we will only give a local answer about the determination of
$p$, working first on a linearized version of the problem, as shown is Section~3. Assuming that $p\in L^{\infty}$ is a given function, we are concerned with the stability around $p$. That is to say $p$ and $u(p)$ are known while $q$ is unknown. 
We can also add that \textit{uniqueness} is a direct consequence of \textit{stability} but historically, 
 uniqueness results were obtain first (see \cite{B-K}) and stability was proved using for instance compactness-uniqueness arguments as in \cite{Y}.\\

In this work we introduce a Carleman weight whose spatial part is similar to the one of the weight function constructed in \cite{B-M-O} for the two-dimensional case. We prove a new Carleman estimate for the Shr\"odinger equation  (see Theorem \ref{Carlemanconvex}) under the hypothesis of {\bf strong convexity} -also called uniform convexity- for the interface (roughly speaking, it means that their curvatures are uniformly bounded from below by a positive constant; see Definition \ref{StrConv}), and some sign for the jump of the main coefficient.  
The following result, based on this Carleman estimate, states the stability of the inverse problem. 

\begin{Thm} \label{ip}
Assume that $\Om_1$ is strongly  convex and $a_1>a_2>0$. Let $\mathcal U$ be a bounded subset of  $L^\infty(\Om)$, $p\in L^\infty(\Om)$ and $r > 0$. If 

$y_0 \in H^1(\Omega)$ 
is real valued (or pure imaginary) and if
$$|y_0(x)| \geq r > 0, ~ a.e. ~ \textit{Êin } \Om,$$
$$y(p) \in H^1(0,T;L^{\infty}(\Om)),$$
then there exists $C = C(\Omega, T, \|p\|_{L^{\infty}(\Om)}, \mathcal U)>0$ such that
$$\|p - q\|_{L^2(\Om)}  \leq C \left\|a_2\dn{y(p)}  - a_2\dn{y(q)}\right\|_{H^1(0,T;L^2(\Gm))}$$
 
$q \in \mathcal U$,
where $y(p)$ is the solution of equation $(1)$ with potential $p$.
\end{Thm}

The main idea is that the {\it nonlinear inverse problem} is reduced to some {\it perturbed inverse problem} which will be solved with the help of a Carleman estimate. In order to obtain such an estimate, we first rewrite (\ref{Anonlin}) as a system of two Schr\"odinger  equations with constant coefficients and solutions $y_1$ and $y_2$, coupled with transmission conditions. 
We then construct a Carleman inequality on each domain with nonzero boundary values on the interface. 
Next, we gather all the terms to construct a global Carleman inequality for the transmission problem. The main point is to look carefully at the interaction of $y_1$ and $y_2$ on the common boundary $\Gamma_1$.

Notice that we state hypothesis for a function which guarantee that it would be a suitable weight function for a Carleman estimate with the only  requirement that the discontinuities of $a$ are located on $\Gamma_1$. 
We shall only construct an explicit weight function for the case of a discontinuous coefficient which 
is constant on each subdomain (i.e. $a_1$ and $a_2$ constants). 
However, we could also construct a weight function for  variable coefficients $a_1(x)$ and $a_2(x)$  such that their traces at the interface are constant,   under additional assumptions of boundedness of $\nabla a_j$ similar to those appearing in \cite{I-Y1} (in order that the corresponding weight function would satisfy  hypothesis $(H_3)$ and $(H_4)$ in section \ref{A2}).\\

This  article is organized as follows. Section~2 is devoted to the proof of an appropriate global Carleman inequality and Section~3 concerns the proof of  the Lipschitz stability of our the inverse problem.


\section{A global Carleman estimate}\label{A2}

In this step, we will show a global Carleman estimate concerning a
function $v=v(x,t)$ equals to zero on $\partial\Omega \times (-T,T)$
and solution of a Schr\"odinger equation with a bounded potential $q=q(x)$. We set the following notations : 
\begin{eqnarray*}
Q= \Om\times (-T,T) && \Om_0 =  \Om_1 \cup \Om_2\\
\Gamma = \partial \Om &&\Gamma_1 = \partial \Om_1 \cap \partial \Om_2 \\
\Sigma = \Gamma\times (-T,T) && \Sigma_1 = \Gamma_1\times (-T,T)
\end{eqnarray*}
and if $u$ is a function defined  in $\Om$, for   $u_j$  we will mean its  restriction to the set $\Om_j$, for each $j = 1,2$.\\

The main hypothesis for the Carleman estimate is the existence of a 
weight function $\psi=\psi(x)$ defined on $\mathbb{R}^N$
such that, on the one hand  it is pseudo-convex with respect to the Schr\"odinger operator
in  each one of the two sub-domains $\Omega_1$ and $\Omega_2$, 
and on the other hand it has a convenient behavior at the interface $\Gamma_1$. \\

Indeed, we will first suppose that $\psi\in C^4(\overline\Om)$ verifies the natural transmission conditions: 
\begin{equation}\label{Tr}             
 \left\{ \begin{array}{cl}
         \psi_1 = \psi_2   &\textnormal{on }\Gamma_1\\
                    a_1 \dns{\psi_1}{1} + a_2 \dns{\psi_2}{2} = 0  &\textnormal{on }\Gamma_1.
          \end{array}
\right. \tag{Tr}
\end{equation}
We will also suppose the following behavior at the interface
\begin{equation} \tag{$H_1$}
\psi(x) = cte \qquad \mbox{  for all  }   x  \in \Gamma_1,
\end{equation}
\begin{equation} \tag{$H_2$}
\dns{\psi_1}1 + \dns{\psi_2}2 < 0 \qquad \mbox{  on  } \Gamma_1.
\end{equation}
In the interior $\Om_0 $ we will need that 
\begin{equation} \tag{$H_3$}
|\nabla \psi| \ge \delta > 0  
\end{equation}
and that $\exists ~\epsilon >0$ such that 
\begin{equation} \tag{$H_4$}
 2 D_a^2\psi(\xi, \bar \xi) + 2 a^2 \la |\nabla \psi \cdot \xi |^2 -a\nabla a\cdot\nabla\psi  |\xi|^2 \ge \epsilon |\xi|^2
\end{equation}
$~\forall \xi \in \CC^n $, where $$D_a^2\psi = \left( a\frac {\partial}{\partial x_i}\left(a\frac {\partial \psi}{\partial x_j}\right) \right)_{1\leq i,j\leq N}.$$


Finally, it will be useful to consider weight functions satisfying $(H_3)$ and $(H_4)$ except  in a neighborhood   of a point.  In this case we will need two weight functions $\psi^1$ and $\psi^2$, each one satisfying $(H_3)$ and $(H_4)$ in $\Om_1 \cup \Om_2 \setminus B_\E(x_1)$  and $\Om_1 \cup \Om_2 \setminus B_\E(x_2)$ respectively, (with $\E$ small enough) 
such that 


\begin{equation} \tag{$H_5$}
\psi^j - \psi^k \ge \delta > 0 \qquad \mbox{  in }  B_\E(x_k)
\end{equation}
for each $ j,k \in \{1,2 \}$ with $j \neq k$.

Summarizing, we set the following

\begin{Def} \label{defwf}
\begin{enumerate} 
\item Let $U \subset \Om$ be an open set such that $\Gamma_1 \subset U$ and let $\psi \in C^4( U \setminus \Gamma_1)$.
 We say that $\psi$  is a transmission weight function for equation $(\ref{Anonlin})$  in $U$  if it  satisfies the conditions \tn{(Tr)}, $(H_1)$ and $(H_2)$ on the interface $\Gamma_1$, and hypothesis  $(H_3)$ and $(H_4)$ in $ U$. \\
\item Let $\psi^1$ and  $\psi^2$ be two functions in $C^4(\Om_1 \cup  \Om_2)$.
 We say that  $(\psi^1, \psi^2 )$ is an $\E$-pair of transmission weight functions for  $(\ref{Anonlin})$ 
if there exist $x_1, x_2 \in \Om_0$ and $\E > 0$ such that 
for each $k  = 1, 2$ the function 
$\psi^k$ is a 
 transmission weight function for  $(\ref{Anonlin})$   in $\Om_0 \setminus \overline{ B_\E(x_k)}$
and  the hypothesis $(H_5)$ is fulfilled.
\end{enumerate}
\end{Def}

Given $\psi$, 
for $s>0$, $\lambda>0$ we define on $Q = \Omega\times (-T,T)$
the following functions:
$$\theta(x,t)=\frac{e^{\lambda\psi(x)}}{(T-t)(T+t)} ~\tn{ and }~
\varphi(x,t)=\frac{\alpha-e^{\lambda \psi(x)}}{(T-t)(T+t)}$$ where
$\alpha>\|e^{\lambda\psi}\|_{L^{\infty}(\Omega)}$.\\

We also define the space
$$Z = \left\{ v\in L^2(-T,T;H^1_0(\Omega))  : 
Lv \in L^2(Q), \dfrac{\partial v}{\partial\nu}\in L^2(\Sigma) \,
\mbox{  and  }  v \mbox{ satisfies  (Tr)} \right\},$$
introduce the following norm in $Z$ 
\begin{equation} \label{norm_s}
\|w\|_{s,\la,\psi}  = s^3 \la^4 \int_{-T}^T\int_\Om \theta^3 |w|^2 dx dt + s \la  \int_{-T}^T\int_\Om \theta |\nabla w|^2 dx dt
\end{equation}
and for $\|\cdot\|_{s,\la,\psi, U}$ we will mean the above terms defined in the set $U\subset\Om$.\\

We finally set
$$Lv = i v' +~\Div (a(x) \nabla u) + qv,$$
$$v=e^{s\varphi} w$$
 and 
$$Pw = e^{-s\varphi}L(e^{s\varphi} w).$$
Hence we have
\begin{eqnarray*}
Pw & = & iw'+is\varphi ' w+~\Div(a \nabla w) + 2s a \nabla\varphi .\nabla w\\
&& + ~sw~\Div(a \nabla \varphi)  +s^2a |\nabla\varphi|^2w + qw\\
& = & P_1w+P_2w +qw
\end{eqnarray*}
where we denoted 
\begin{eqnarray*}
&&P_1w=iw'+~\Div (a\nabla  w) +s^2a|\nabla\varphi|^2w,\\
&&P_2w=is\varphi 'w+2sa \nabla\varphi .\nabla w+s~\Div(a \nabla \varphi) w.
\end{eqnarray*}


Our main result is the following
\begin{Thm} \label{Carleman} 
Suppose  there exists for some $\E > 0$ an $\E$-pair of  transmission weight functions   $(\psi^1, \psi^2)$ belonging to $C^3(\Om_1 \cup \Om_2)$. Let  $\theta^k$, $\varphi^k$ and $w^k$ be the corresponding functions defined for $\psi^k$ as we did before. We also define
$$\Sigma_+^{\psi^k} = \left\{ (x,t) \in \Gamma\times(-T,T) \, : \nabla \psi^k(x,t) \cdot \nu(x) > 0 \right\}$$
Then there exists $C>0$, $s_0>0$ and $\la_0 >0$ such that
 \begin{gather}
 \sum_{k=1}^{2} \left(  
 \left\|P_1^{\psi^k}(w^k) \right \|^2_{L^2(Q)} \right.  + \left. \left \|P_2^{\psi^k}(w^k) \right \|^2_{L^2(Q)}
+   \left\|w^k \right\|_{_{\la, s, \psi^k}}^2   \right)  \nonumber  \\
 \leq  C~ \sum_{k=1}^2 \left(   \left \|P^{\psi^k}(w^k) \right\|^2_{L^2(Q)}   +  s \la  \iint_{\Sigma_+^{\psi^k}} \theta^k \left|a \dn{w^k} \right|^2  \right)  \label{ineqCarlem}
 \end{gather}
for all $v \in Z$,  $\la \geq \la_0$ and $s \geq s_0$. 
\end{Thm}



\subsection{Formal computations}

We have
\begin{eqnarray*}
\iint_{Q}|Pw - qw|^2\,dxdt&=&
\iint_{Q}|P_1w|^2\,dxdt
+\iint_{Q}|P_2w|^2\,dxdt\\
&&+2Re \iint_{Q}P_1w\overline{P_2w}\,dxdt,
\end{eqnarray*}
where $\overline{z}$ is the conjugate of $z$ and $Re(z)$ its real part.\\
As $v\in L^2(-T,T;H^1_0(\Omega))$ and $v'\in L^2(-T,T;H^{-1}(\Omega))$
(because $Lv \in L^2(Q)$), we have $v\in C([-T,T];L^2(\Omega))$ and $w\in C([-T,T];L^2(\Omega))$ with $w(x,\pm T)=0$.\\


We  look for lower  bounds for $$Re\iint_{Q}P_{1}w\overline{P_2w}\,dxdt = \left<P_1w,P_2w\right>_{L^2}$$ 

We set  $\left<P_1w,P_2w\right>_{L^2} = \sum\limits_{i,j=1}^3 I_{i,j}$, where  $I_{i,j}$ is the integral of the product of the $i$th-term of $P_1w$ and the $j$th-term of $P_2w$. 
The properties of $w$ and some integrations by parts allow to write the following equalities.

To begin with, we have

\begin{eqnarray*}
&&I_{11}=Re\iint_{Q}iw'(-is\varphi '\overline{w})\,dxdt=
-\frac s2\iint_{Q}\varphi''|w|^2\,dxdt.
\end{eqnarray*}
Applying the identity $Im(z) = - Im(\overline{z})$ for $z=2s\lambda~\ds\int_{-T}^{T}
\ds\int_{\Omega}\theta a \nabla\psi \cdot \nabla\overline{w}w'\,dxdt $ we obtain:
\begin{eqnarray*}
I_{12}&=&Re\iint_{Q}iw'(2s a \nabla\varphi \cdot \nabla\overline{w})\,
dxdt\\
&=&s\lambda~ Im\iint_{Q}\theta(~\Div (a \nabla\psi)+\lambda a |\nabla\psi|^2)w\overline{w'}\,dxdt \\ 
&&-~ s\lambda~ Im\iint_{Q} a \theta'w
\nabla\psi \cdot \nabla\overline{w}\,dxdt\\
&& + s \la \iint_\Sg \overline w a \theta w' \dn{\psi}\,d\sigma dt
\end{eqnarray*}

We also have
\begin{eqnarray*}
I_{13}&=&Re\iint_{Q}iw's~\Div (a\nabla \varphi ) \overline{w} \,dxdt\\
&=& -s\lambda~Im\iint_{Q}\theta(~\Div(a\nabla\psi)+\lambda a |\nabla\psi|^2) w\overline{w}'dxdt,
\end{eqnarray*}
\begin{eqnarray*}
I_{21}&=&Re\iint_{Q}~\Div(a \nabla w)(-is\varphi '\overline{w})\,dxdt\\
& =& s\lambda~Im\iint_{Q} a \theta'\overline{w}
\nabla\psi \cdot \nabla w\,dxdt 
+  s~Im\iint_{\Sigma}\varphi'\overline{w}a\frac{\partial w}{\partial \nu}\,d\sigma dt,
\end{eqnarray*}
and
\begin{eqnarray*}
I_{22}&=&Re\iint_{Q}~\Div(a \nabla w)(2s a \nabla\varphi \cdot \nabla \overline{w}) \,dxdt\\
&=&-~s\lambda\iint_{Q}\theta a |\nabla w|^2 ( ~\Div(a \nabla \psi) + \lambda  a |\nabla \psi|^2)\,dxdt \\
&&-~s\lambda\iint_{Q}\theta a |\nabla w|^2\nabla a \cdot \nabla \psi\,dxdt  +~2s\lambda^2\iint_{Q}\theta a^2 |\nabla\psi.\nabla w|^2\,dxdt \\
&&+~2s\lambda~Re\iint_{Q}\theta D_a^2(\psi)(\nabla w,\nabla \overline w) \,dxdt\\
&&-~2 s\lambda\iint_{\Sigma}\theta a^2
(\nabla \psi \cdot \nabla \overline w) \frac{\partial w}{\partial \nu} \,d\sigma dt
+  s\lambda\iint_{\Sigma}\theta a^2
|\nabla w|^2 \frac{\partial \psi}{\partial \nu} \,d\sigma dt.
\end{eqnarray*}
Using integrations by parts we obtain
\begin{eqnarray*}
I_{23}&=&Re\iint_{Q}~\Div(a \nabla w)
s~\Div(a \nabla \varphi)\overline{w}\,dxdt\\
&=&s\lambda \iint_{Q} |\nabla w|^2 \theta(a ~\Div(a \nabla \psi) 
+\lambda |a \nabla\psi|^2)
\,dxdt \\
&&-~\frac{s\lambda}2 \iint_{Q} |w|^2 ~\Div(  \theta a \nabla ~\Div(a \nabla \psi))\,dxdt   \\
&& + ~  \frac{s\lambda}2 \iint_\Sigma |w|^2\theta a \nabla ~\Div(a\nabla \psi)\cdot \nu \,d\sigma dt\\
&&-~\frac{s\lambda^2}2 \iint_{Q} |w|^2(  ~\Div(a \theta \nabla (a|\nabla\psi|^2)) + ~\Div(a \theta ~\Div(a\nabla \psi) \nabla \psi) )
\,dxdt \\
&& + \frac{s\lambda^2}2 \iint_\Sigma |w|^2  \theta a( ~\Div(a  \nabla \psi) \nabla \psi + \nabla(a|\nabla \psi|^2) ) \cdot \nu\,d\sigma dt\\
&&-~\frac{s\lambda^3}2 \iint_{Q} |w|^2 ~\Div( a^2 \theta|\nabla\psi|^2\nabla \psi)\,dxdt
+ \frac{s\lambda^3}2 \iint_\Sigma |w|^2 a^2 \theta | \nabla\psi|^2 \dn{\psi}\,d\sigma dt\\
&&- \la s Re \iint_\Sg \overline w \theta a (~\Div(a \nabla \psi) + \la |\nabla \psi|^2a) \dn{w}\,d\sigma dt .
\end{eqnarray*}
and we obviously have
\begin{eqnarray*}
I_{31}&=&Re\iint_{Q}s^2a |\nabla\varphi|^2w(-is\varphi ' \overline w)
\,dxdt=0,
\end{eqnarray*}
\begin{eqnarray*}
I_{32}&=&Re\iint_{Q}s^2a |\nabla\varphi|^2w(2s a \nabla\varphi \cdot \nabla \overline w)\,dxdt\\
&=&s^3\lambda ^3~\iint_{Q} |w|^2 \theta^3a \left( |\nabla\psi|^2  ~\Div(a\nabla\psi)+ 2 a D^2(\psi) (\nabla \psi, \nabla \psi )\right) \,dxdt  \\
&&+~s^3\lambda ^3~\iint_{Q} |w|^2 \theta^3a  |\nabla \psi |^2 \nabla a \cdot \nabla \psi  \,dxdt\\
&&+~3s^3\lambda ^4~\iint_{Q} |w|^2 \theta^3 a^2 |\nabla\psi|^4\,dxdt
- \lambda^3 s^3 \iint_\Sigma |w|^2a^2|\nabla \psi|^2 \theta^3 \dn{\psi}\,d\sigma dt,
\end{eqnarray*}
and
\begin{eqnarray*}
I_{33}&=&Re\iint_{Q}s^2a |\nabla\varphi|^2w
(s~\Div(a \nabla\varphi ) \overline w)\,dxdt\\
&=&-~s^3\lambda ^3~\iint_{Q}|w|^2 \theta^3 a |\nabla\psi|^2~\Div(a \nabla\psi) \,dxdt\\
&&-~s^3\lambda ^4~\iint_{Q}|w|^2 \theta^3a^2 |\nabla\psi|^4\,dxdt.
\end{eqnarray*}
Then we have 
$~~Re \displaystyle\iint_{Q}P_{1}w\overline{P_{2}w}\,dxdt = F(w) + G(\nabla w) + J + X_1~~$
where we define
\begin{eqnarray}
 F(w) &=& 2 s^3 \la^4\iint_Q |w|^2 \theta^3 a^2 |\nabla \psi|^4\,dxdt, \label{DefF} \\
G(\nabla w) 
& = &  2 s \la^2 \iint_Q  \theta a^2 |\nabla \psi \cdot \nabla w|^2 \,dxdt  + 2 s \la Re \iint_Q \theta  D_a^2\psi(\nabla w, \nabla \overline w)\,dxdt \nonumber \\
& & -   s \la \iint_Q |\nabla w|^2 \theta a \nabla a \cdot \nabla \psi\,dxdt, \label{DefG}
\end{eqnarray}  
$J$ as the sum of all boundary integrals
\begin{eqnarray*} 
J & = & s \la Im \iint_\Sg a \theta w' \overline w \dn{\psi}\,d\sigma dt \\
& & s  Im \iint_\Sg \varphi'\overline w a \dn{w} \,d\sigma dt  \\
& & - 2 s \la Re \iint_\Sg a^2 \theta \nabla \psi \cdot \nabla \overline w \dn{w}\,d\sigma dt  \\
& & + s \la  \iint_\Sg a^2 \theta |\nabla w|^2 \dn{\psi}\,d\sigma dt  \label{listJ}  \\
& & -  s \la Re \iint_\Sg \overline w ~\Div(a \theta \nabla \psi)a \dn{w}\,d\sigma dt  \\
& & - s^3 \la^3  \iint_\Sg a^2 \theta^3 |w|^2 |\nabla \psi|^2 \dn{\psi}\,d\sigma dt \\
& & + \frac{s \la^3}2  \iint_\Sg  a^2 \theta |w|^2 |\nabla \psi|^2 \dn{\psi}\,d\sigma dt  \\
& & + \frac{s \la^2}2 \iint_\Sg a^2 \theta |w|^2  \nabla( |\nabla \psi|^2 )\cdot \nu \,d\sigma dt \\
& &  + \frac{s \la^2}2  \iint_\Sg a \theta |w|^2 ~\Div(a \nabla \psi)\dn{\psi} \,d\sigma dt\\
& & +  \frac{s \la}2 \iint_\Sg a \theta |w|^2 \nabla(~\Div(a \nabla \psi))\cdot \nu\,d\sigma dt
\end{eqnarray*}
and $X_1$  as  the  sum of  all the remaining integrals in $\Omega$.\\

Moreover, if $U \subset \Om$ is an open set, we will write $F_U(w)$ to denote  the sum of integrals from the  definition of  $F(w)$  taken in the set $U$, and the same for $G$, $X_1$...  \\

Noticing that \begin{eqnarray*}
&\bullet& 2s\lambda ~ Im\iint_{Q}\theta' a w \nabla\psi.\nabla\overline{w}\,dxdt \\
&&\quad\leq s\lambda \iint_{Q}a^2(\theta')^{\frac {1}{2}} |\nabla\psi.\nabla w|^2\,dxdt
+ ~s\lambda \iint_{Q}(\theta')^{\frac {3}{2}}|w|^2\,dxdt,\\ 
& \bullet& a\in W^{2,\infty} (\Omega) \tn{ and } \psi\in C^4(\overline\Om)\\ 
& \bullet& |\theta|\leq C\theta^3,~|\theta'|\leq C\theta^2 \tn{ and} |\varphi''|\leq C\theta^3 \tn{ on } (-T,T)\times\Omega,~C=C(T)>0.
\end{eqnarray*}
it is then easy to prove, from simple calculations, that the ``negligible" terms $X_1$ indeed satisfy
\begin{equation} \label{mayX1}
 |X_1| \le C s\la \iint_Q a^2 \theta |\nabla \psi \cdot \nabla w|^2
+ C s \la^4 \iint_Q \theta |w|^2 
+  C s^3 \la^3 \iint_Q \theta^3|w|^2.
\end{equation}

\subsection{Proof of the Carleman estimate}

In this part of the paper, we prove Theorem~\ref{Carleman}. We apply the above computations in each one of the domains $\Om_1$ and $\Om_2$ and we sum up all the terms. Since the interface  $\Gamma_1$ has null $\RR^N$-measure, we get an estimate in all the set $\Om$, plus the boundary terms from $\partial \Om$, and from the interface itself, where appear terms coming from both $\Om_1$ and $\Om_2$. \\
Given the hypothesis we have assumed, we prove in the following propositions that we can deal with all this terms. 
In the sequel, $C$ denotes a generic constant, depending on $T$ and $\Omega$.

\subsubsection{The interior}

Recall the norm $\|\cdot\|_{s,\la,\psi}$ defined in  $(\ref{norm_s})$ and $F(w)$, $G(\nabla w)$ defined in
$(\ref{DefF})$, $(\ref{DefG})$.


\begin{Prop} \label{propint}
Suppose that $U \subset \Om$ is a open set and $\psi$ satisfies $(H_3)$ and $(H_4)$ in $\Om_1 \cup \Om_2 \setminus U$. Then there exist $\gamma > 0$ , $C \in \RR$, $ s_0$ and $\la_0$  such that for all $v \in Z$,
$$ F(w) + G(\nabla w) + X_1 \ge \gamma \|w\|_{s,\la,\psi}  - C \| w \|_{s,\la,\psi,U}  $$
$\forall s \ge s_0$ and $\forall \la \ge \la_0$.
\end{Prop}

\dem 

First, merely by the fact  that $\psi \in C^4(\overline\Om)$, we have
\begin{equation} \label{mayU}
|F_U(w)| + |G_U(\nabla w)| + |X_{1,U}| \le C \|w\|_{s, \la, \psi, U} 
\end{equation}
for all $v \in Z$.

Now, from $(\ref{mayX1})$ and $\psi$ satisfying $(H_3)$ we get that for $s$ and $\la$ large enough,
\begin{equation*}  \label{mayX12}
|X_1| \le  s \la^2 \iint_Q \theta a^2 |\nabla \psi \cdot \nabla w|^2 
             +  s^3\la^4 \iint_Q |w|^2\theta^3 a^2 |\nabla \psi|^4 
\end{equation*}
Hence, If $\psi$ satisfies $(H_3)$ and $(H_4)$ in $\Om_* = \Om_1 \cup \Om_2 \setminus U$ 
we get that $\forall\la \ge \la_0$, $s \ge s_0$ and $v \in Z$,
\begin{eqnarray}
  F_{\Om_*}(w)  & + &  G_{\Om_*}(\nabla w) + X_{1,\Om_*} \nonumber  \\
 &  \ge  &  F_{\Om_*}(w) + G_{\Om_*}(\nabla w) - |X_{1,\Om_*}|  \nonumber \\
&  \ge  &   F_{\Om_*}(w) + G_{\Om_*}(\nabla w) - s \la^2 \int_{-T}^T\int_{\Om_*} \theta a^2 |\nabla \psi \cdot \nabla w|^2 \nonumber\\
&&- ~s^3\la^4 \int_{-T}^T\int_{\Om^*} |w|^2\theta^3 a^2 |\nabla \psi|^4 \nonumber \\
&  \ge  &    s^3\la^4 \int_{-T}^T\int_{\Om_*} |w|^2\theta^3 a^2 |\nabla \psi|^4 
+ s \la^2 \int_{-T}^T\int_{\Om_*} \theta a^2 |\nabla \psi \cdot \nabla w|^2 \nonumber \\
& & + ~2 s \la \mbox{Re} \int_{-T}^T\int_{\Om_*} \theta D_a^2\psi(\nabla w, \nabla \bar w) -  s \la \int_{-T}^T\int_{\Omega_*} |\nabla w|^2 \theta a \nabla a \cdot \nabla \psi   \nonumber  \\
&  \ge  &   s^3\la^4 \int_{-T}^T\int_{\Om_*} |w|^2\theta^3 a^2 |\nabla \psi|^4 +  \epsilon s \la \int_{-T}^T\int_{\Om_*} \theta   |\nabla w|^2   \nonumber \\
& \ge & \gamma \|w\|_{s, \la, \psi, {\Om_*}}.  \label{min*}
\end{eqnarray}

From (\ref{mayU})  and (\ref{min*})  we get the desired result.
\QED

\subsubsection{The boundary}

By definition  we have $w = 0 $ on the exterior boundary $\Sg$ for each  $v \in Z$. 
Therefore, $\nabla w|_{\Sg} = \dn{w} \nu$ 
and if we choose the legitimate notation $J=  J_{\Sigma} + J_{\Sigma_1} $, we get here
\begin{eqnarray} 
 J_{\Sigma}   & = &- 2 s \la Re \iint_\Sg a^2 \theta \nabla \psi \cdot \nabla \overline w \dn{w}\,d\sigma dt 
 + s \la  \iint_\Sg a^2 \theta |\nabla w|^2 \dn{\psi}\,d\sigma dt \nonumber\\
& = &-s \la \iint_{\Sg} \theta  \left|a\dn{w} \right|^2 \dn{\psi}\,d\sigma dt  \nonumber \\
& \geq & -s \la \iint_{\Sigma_+} \theta   \left|a\dn{w} \right|^2\dn{\psi}\,d\sigma dt \nonumber \\
& \geq & -s \la \left\|\dn{\psi}\right\|_{L^{\infty}(\Sg)} \iint_{\Sigma_+}  \theta   \left|a\dn{w} \right|^2\,d\sigma dt \nonumber \\
& \geq &  -s \la C\iint_{\Sigma_+}  \theta   \left|a\dn{w} \right|^2\,d\sigma dt\label{bd}
 \end{eqnarray}					   
where we have denoted 
$\Sigma_+ = \left\{ (x,t) \in \Gamma\times(-T,T) \, : \nabla \psi(x,t) \cdot \nu(x) > 0 \right\}$.

\subsubsection{The interface}

We compute the sum of the integrals  on the interface $\Sg_1$, 
 We write $\ds{J_{\Sg_1} = \sum_{k=1}^{10}J_k}$,  enumerating the terms in the same order of the list in $(\ref{listJ})$. For each $k = 1, \ldots, 10$ we denote as $[J_k]$ the sum of the  $k$-term coming from the integrations by parts in $\Omega_1$ with the corresponding one from $\Omega_2$.

\begin{Prop} \label{propGamma1}
If $\psi$ satisfies  hypothesis $(H_1)$, $(H_2)$ and \tn{(Tr)}  then  there exist $\la_0$ and $s_0$ such that
\begin{equation}
J_{\Sg_1} = \sum_{k=1}^{10} (J_k(w_1) + J_k(w_2) ) \ge 0
\end{equation}
for all $v \in Z$, $\forall \la \ge \la_0$, $s \ge s_0$.
\end{Prop}

\dem

It is not difficult to check that $[J_k] = 0$ for $k=1,2$ since $\psi$ and $w$ satisfy the transmission conditions (Tr).
Moreover, $\psi$ is constant on the interface and then we obtain  
$\ds{\nabla \psi \cdot \nabla \overline w = \dn{\psi}\dn{\overline w}}$ on $\Gamma_1$. 
Therefore, thanks to $(H_2)$ we get
\begin{eqnarray*}
[J_3] & = & - 2 s \la  \iint_{\Sg_1}  \theta \left| a_1 \dns{w_1}1\right|^2 \left(  \dns{\psi_1}1 + \dns{\psi_2}2 \right) \,d\sigma dt  \\
 & \ge &  s \la \delta  \iint_{\Sg_1}   \theta \left| a_1 \dns{w_1}1\right|^2\,d\sigma dt .
\end{eqnarray*}

By mean of the orthogonal  decomposition
$\ds{ \nabla w = \dn{w} \nu+ \nabla_\tau w}$, where $\nabla_\tau w$ is the projection of  $\nabla w$ 
on the tangent hyper-plane of $\partial \Omega_1$, and from 
hypothesis $(H_2)$ and $(H_3)$ and the fact that  
 $ \nabla_\tau w_1 =  \nabla_\tau w_2 $
we obtain
\begin{eqnarray*}
[J_4] & = & s \la  \iint_{\Sg_1} \theta \left|a_1 \dns{ w_1}1\right|^2  \left(   \dns{\psi_1}1 +   \dns{\psi_2}2   \right)\,d\sigma dt  \\
& & + ~s \la  \iint_{\Sg_1} \theta \left|\nabla_\tau w\right|^2 \left( a_1^2 \dns{ \psi_1}1 + a_2^2 \dns{ \psi_2}2   \right)\,d\sigma dt   \\
& =  &  s \la  \iint_{\Sg_1} \theta \left|a_1 \dns{ w_1}1\right|^2  \left(   \dns{\psi_1}1 +   \dns{\psi_2}2   \right)\,d\sigma dt  \\
& & + ~s \la  \iint_{\Sg_1} \theta \left|\nabla_\tau w \right|^2 a_1a_2 \left( - \dns{ \psi_1}1 - \dns{ \psi_2}2   \right)  \,d\sigma dt \\
& \ge  &  s \la  \iint_{\Sg_1} \theta \left|a_1 \dns{ w_1}1\right|^2  \left(   \dns{\psi_1}1 +   \dns{\psi_2}2   \right)\,d\sigma dt  \\
 & \ge & 2 s \la \delta  \iint_{\Sg_1}  \theta \left| a_1 \dns{w_1}1\right|^2 \,d\sigma dt 
 \end{eqnarray*}

\begin{eqnarray*}
[J_6] & = &
 - s^3 \la^3 \iint_{\Sg_1}  \theta^3 |w_1|^2 \left|a_1\dns{\psi_1}1\right|^2 \left( \dns{\psi_1}1 + \dns{\psi_2}2\right) \,d\sigma dt \\
&  \ge & s^3 \la^3 \delta^2 \iint_{\Sg_1}  \theta^3 a_1^2 |w_1|^2 \,d\sigma dt 
\end{eqnarray*}

Since $a\in W^{2,\infty} (\Omega)$ and  $\varphi\in C^4(\overline\Om)$, we also have
\begin{eqnarray*}
\left| [J_5] \right|& \le & C s^2 \la^3 \iint_{\Sg_1}|w_1|^2\theta^3 \,d\sigma dt  +  C \la \iint_{\Sg_1}\left|a\dns{w_1}1\right|^2\theta \,d\sigma dt 
\end{eqnarray*}
and
\begin{eqnarray*}
\left| \sum_{k=7}^{10} [J_k] \right| & \le & C s^2 \la^3 \iint_{\Sg_1}|w_1|^2\theta^3 \,d\sigma dt . 
\end{eqnarray*}

Thus, for $s$ large enough, we get the desired result
$$\sum_{k=4}^{10} [J_k]  \ge  
(s\delta - C) \left(  s^2 \la^3 \iint_{\Sg_1}|w_1|^2\theta^3 +   \la \iint_{\Sg_1}|a\dns{w_1}1|^2\theta \right) \ge  0$$
\QED

 \subsubsection{Carrying all together.} 
 From $(\ref{bd})$ and Propositions \ref{propint} and \ref{propGamma1}  we obtain 
 \begin{equation}
 \left\|w\right\|_{_{s, \la,\psi}}^2   - C \left\|w \right\|_{_{s, \la,\psi, U}}^2 
  - s \la C \iint_{\Sigma_+} \theta \left |a\dn{w} \right |^2    \leq   C Re \left< P_1(w), P_2(w) \right>_{L^2}  \label{ineg1}
\end{equation} 

Adding $ \dfrac{C}{2} \left( |P_1(w)|^2_{L^2} + |P_2(w)|^2_{L^2}\right) $ to 
both sides of $(\ref{ineg1})$ we obtain
\begin{multline*}
\dfrac{C}{2}\left(  \left|P_1(w) \right |^2_{L^2} + \left|P_2(w) \right|^2_{L^2}\right)
 +   \left\|w \right\|_{_{s, \la,\psi}}^2 \\
  - C  \left\|w \right\|_{_{s, \la,\psi,U}}^2 
      - s \la C \iint_{\Sigma_+} \theta \left |a\dn{w} \right |^2   
         \leq  C \left |P(w) -  q w\right |^2_{L^2},  
 \end{multline*} 
what means that for all $s \geq s_2$  and $\la \geq \la_2$, since $C>0$ is a generic constant,
\begin{multline}
  |P_1(w)|^2_{L^2} + |P_2(w)|^2_{L^2}
+ \left \|w \right\|_{_{s, \la,\psi}}^2 \\
   \leq  C |P(w) -  q w |^2_{L^2}  
 + C \left\|w\right\|_{_{s, \la,\psi,U}}^2 
  +    s \la C  \iint_{\Sigma_+} \theta \left|a\dn{w}\right |^2.  \label{Carl0}
\end{multline}

Now, if $(\psi^1, \psi^2)$ is an $\E$-pair of  transmission weight functions 
(see Definition~\ref{defwf}), we have an estimate like $(\ref{Carl0})$ for each $\psi^k$ with $U =  B_\E(x_k)$ where $x_j \in \Om$, $j=1,2$ and $\E > 0$. 

We sum up both estimates and we can show that the left hand side of each inequality can absorb the right hand side term $\| \cdot \|_{s,\la,\psi^k,B_\E(x_k)}$ from the other inequality provided that $\E$ is small and $\la$ is large enough.
Indeed, by assumption we have that  $ \psi^2 -  \psi^1 > \delta > 0$ in $\Bi$. \\
Then, by taking $\la$ large enough we have
$$e^{\la( \psi^2 -  \psi^1)} > 2C \hspace{1cm}  \mbox{in}  \quad \Bi$$
i.e.
$$C  \theta^1 < \frac{1}{2}  \theta^2  \hspace{1cm}  \mbox{in} \quad \Bi$$
and we conclude that  $\|w^1 \|_{\psi^1,\Bi}$ of the right and side is absorbed by the term $ \|w^2 \|_{_{\psi^2}}$ of the left hand side.
It is clear that an analogous result is true by interchanging $\psi^1$ and $\psi^2$. Theorem~\ref{Carleman} is proved.
\QED

\subsection{Particular case.}

In this part of the work, we construct explicit weight functions adapted to particular discontinuous coefficients. 

We need the following definition.

\begin{Def}\label{StrConv}
We say that the  open, bounded and convex set $U \subset \RR^N$ ($N \geq 2$) 
 is {\bf strongly convex} if $\partial U$ is of class
 $C^2$ and all the principal curvatures are strictly positive functions on $\partial U$.
\end{Def}

\begin{rem}  \label{planes}
Let us note that $U \subset \RR^N$ is strongly convex  if and only if 
for all plane $\Pi \subset \RR^N$ intersecting $U$,  the curve   $\Pi \cap \partial U$  has strictly positive curvature at each point. 
In particular, a strongly convex set is geometrically strictly convex.
\end{rem}

\medskip

We assume that $\Om_1\subset\Omega$ is a strongly convex domain  with boundary $\Gamma_1$ of class $C^3$, 
and we set $\Om_2 = \Om \setminus \overline{\Om}_1$. 
Thus, we have $\partial \Om_2 = \Gm \cup \Gm_1$, where this is a disjoint union.
We deal with the case where $a$ is locally constant
\begin{equation} \label{a1a2}
a(x) = \left\{ \begin{array}{cc}
a_1    &x\in \Om_1 \\
a_2   &x\in \Om_2
\end{array}
\right. 
\end{equation}
with $a_j  > 0$ for $j = 1,2$. 

In order to construct a convenient weight function, take $x_0 \in \Om_1$ and for each $x\in\Om\setminus \{x_0\}$ define $\ell(x_0,x) = \{ x_0 + \la(x - x_0) \, : \,  \la \geq 0 \}$.  
Since $\Om_1$ is convex there is exactly one point $y(x)$ such that
$y(x) \in \Gm_1 \cap \ell(x_0,x)$.
Thus, we can define the function 
$\fun{\rho}{\Om\setminus \{x_0\}}{\RR}^+$ by:
\begin{equation}\label{defrho}
\rho(x) = |x_0 - y(x)|. 
\end{equation}
Let $\E > 0$ be such that $\overline B_\E\subset\Omega_1$ (and small enough in a sense we will precise later) and let $0<\E_1<\E_2<\E$. Then we consider a cut-off function 
$\eta \in C^{\infty}(\RR^N)$ such that 
$$0 \leq \eta \leq 1,\qquad
\eta = 0 \textrm{ in } B_{\E_1}(x_0),\qquad
\eta = 1 \textrm{ in } \Omega\setminus \overline B_{\E_2}(x_0).$$
For each $j \in \{ 1,2\}$ we take $k$ such that $\{j,k \} = \{1,2\}$ and we define the following
functions in the whole domain $\Omega$ 
$$\psi_j(x)  =  \eta(x) \frac{a_k}{\rho(x)^2}|x - x_0|^2 + M_j  
\qquad  x  \in \Om,$$
where $M_1$ and $M_2$ are positive numbers such  that 
\begin{equation} \label{aM}
a_1 -  a_2 =  M_1 -  M_2.
\end{equation}
Then, the weight function we will use in this work is
\begin{equation}\label{phi}
\hspace{1cm} 
\psi(x) = 
\left\{ \begin{array}{lrr}
 \psi_1(x)     
 & x \in \Om_1  & \vspace{0,3cm}\\ 
\psi_2(x)    
& x \in \Om_2. &
\end{array}
\right.
\end{equation}

Throughout the paper, we will use the notations $~\bar a(x) =  a_2 \mathbf{ 1}_{\Om_1}(x) + a_1 \mathbf{1}_{\Om_2}(x)~$and $~M =  M_1 $ on $\Om_1$ and $M= M_2 $ on $\Om_2$, so that we can write 
$$\ds{\psi(x) = \eta(x) \bar a (x)\frac{|x - x_0|^2}{\rho(x)^2} + M}.$$

As we can see in the following result, the main property of the weight function is a consequence of the strong convexity of the interior domain $\Om_1$.
\begin{Lem}\label{positDef}
If $\Om_1 \subset \RR^N$ is strongly convex and if the function $\ds{\fun{\mu}{\RR^N}{\RR^+} }$
is defined by $\ds{\mu(x) = \frac{|x-x_0|}{\rho(x)}}$ then 
$D^2\mu^2(x)$ is positive definite for all $x \in \RR^N\backslash \{ x_0\}$, 
uniformly in bounded subsets of $\RR^N\backslash \{ x_0\}$.
\end{Lem}

\dem

We shall deduce this Lemma from  well-known  properties of  compact convex subsets of $\RR^N$   (called {\it convex bodies}). However, for the sake of completeness of this paper, we include  in the Appendix a self-contained   proof of this result.

Assuming without lost of generality that $x_0  = 0$, it is not difficult to see that $\mu$ is the gauge function of the convex set $\overline{\Om_1}$ (in other words, $\mu$ is a seminorm whose unit ball is $\Om_1$, see \cite{Schn}, p. 43, and section 2.3 of \cite{CraMa}). 

The proof that  $\mu$ is a convex function of class $C^2$ (hence $D^2\mu \ge 0$) can be read in  \cite{CraMa} (Theorem 2.1).
Moreover, it is proved that  for each $x \in \RR^N \backslash \{0\}$ the only null   eigenvalue of $D^2\mu(x)$ corresponds to the direction  $x$ (which is the radial direction). The others eigenvalues, as functions of $x$, are
bounded below by a positive constant, uniformly in   $x \neq 0$ given in a bounded subset of $\RR^N$.

 Thus,  there exists $\delta > 0$ such that for all $x \in \overline \Omega$ we have
\begin{equation} \label{ortx}
D^2\mu(x)(v,v) \ge \delta |v|^2 \, \, \, \forall  \, v \in x^\bot = \{ y \in \RR^n \, : \, y \cdot x = 0  \}.
\end{equation}  

On the other hand, we have 
$ \nabla \mu(x) =  \ds\frac1{\rho(x)} \frac{x}{|x|} + |x| \nabla \ds\frac1\rho (x).$
Since $\rho$ is constant in the radial direction, we get $x \cdot  \nabla \ds\frac{1}{\rho}(x) = 0$. Hence we deduce that
\begin{equation}\label{rad}
x \cdot \nabla \mu(x) =  \frac{|x|}{\rho(x)} = \mu(x) \neq 0.
\end{equation}

Take  $x , v \in \RR^N \backslash \{0\}$. Then $v = v_1 \frac{x}{|x|} + v_2 y$, where $y$ is an unitary element of $x^\bot$.   In view of the fact that
$D^2\mu^2(v,v) = 2 \mu D^2\mu(v,v) + 2 |v \cdot \nabla \mu|^2$, from $(\ref{ortx})$, $(\ref{rad})$ and $D^2\mu(x)(x,x) = 0$ we get
\begin{eqnarray*}
D^2\mu^2(x)(v,v) & \ge & 2\mu(x) D^2\mu(x)(v_2y,v_2y)  +  2\left| v_1\frac{x}{|x|} \cdot \nabla \mu(x) \right|^2 \\
  			    & \ge & 2\mu(x) \delta v_2^2  +  2\frac{v_1^2}{|x|^2} \mu^2(x) \\  
  			    & =  & 2\mu(x) \delta v_2^2  +  \frac{2}{\rho(x)^2} v_1^2 \\  
     			    & \ge & \delta_1 ( v_1^2  + v_2^2) \\  
\end{eqnarray*}
and we conclude that $D^2\mu^2(x)$ is positive definite. 

\QED

\medskip 

Assuming the additional hypothesis about the sign of the  jump on the interface, we can prove that the functions we have defined work as a weight function: 
\begin{Prop}\label{psi12} Let $\Om_1$ be an open and bounded set in $\RR^N$ with smooth boundary,  and $ a_1$, 
$ a_2 $  real numbers such that:
\begin{enumerate}
\item $\Om_1$ is strongly convex.
\item $0 < a_2 < a_1$.
\end{enumerate}
Then, for each pair of points $x_1, x_2 \in \Om_1$, there exists $\E > 0$ such that the above construction gives up an $\E$-pair of transmission weight functions $(\psi^1,\psi^2)$ in the sense of Definition \ref{defwf}.
\end{Prop}

\dem

For $x_0 \in \Om_1$ let $\psi$ be the function constructed as above and defined by (\ref{phi}). If $x \in \Gamma_1$ we have $\rho(x) = |x-x_0|$ and $\psi_j(x) = a_k + M_j$. From $(\ref{aM})$ we get $\psi_1 = \psi_2 = c $ on $\Gamma_1$. 
Moreover, if $x \in \Gamma_1$ we have 
$$a_1 \nabla \psi_1(x) = a_1 a_2 \nabla \left( \frac{|x-x_0|^2}{\rho(x)^2} \right) =  a_2 \nabla \psi_2(x) .$$ Hence (Tr) and ($H_1$) are satisfied (recall that $\nu_1 = -\nu_2$ on $\Gamma_1$). \\

On the other hand,   since $\Gamma_1 $ is a level set  of $\psi_1$, then 
$\psi_1(x) < a_2 + M_1 < \psi_1(y)$ for any $x \in \Om_1$ and $y \in \Om_2$, and we have $\dns{\psi_1}{1} > 0$ on $\Sigma_1$ and
$$ \dns{\psi_1}{1} + \dns{\psi_2}{2} = \dns{\psi_1}{1}  \left( 1 - \frac{a_1}{a_2} \right) < 0$$
what gives ($H_2$).\\

For $x\in \Om_0 \setminus B_\E(x_0)$, denoting $ c(x) = \frac{\bar a}{\rho^2(x)} $,
 we get $\nabla \psi = 2c(x)(x-x_0) + |x-x_0|^2\nabla c(x)$.  
By construction $c(x)$ is constant in the direction of $x-x_0$, hence 
$$(x-x_0) \cdot \nabla c(x) = 0 $$
and then
\begin{eqnarray*}
|\nabla \psi|^2 & = & 4c^2(x)|x-x_0|^2 + |x-x_0|^4|\nabla c(x)|^2 \\
		&  \geq &  4c^2(x)|x-x_0|^2.
\end{eqnarray*}
Thus we have 
$$|\nabla \psi|^2   \geq   4 \left( \frac{\bar{a}}{\textrm{diam}(\Omega)^2} \right)^2 \E^2.$$
in $\Om_0 \setminus B_\E(x_0)$ and  $\psi$ satisfies $(H_3)$ in that set. \\

Property $(H_4)$ is deduced from Lemma \ref{positDef}.\\

One can notice that $x_0$ can be arbitrarily chosen in $\Om_1$ since it is  convex.
Therefore, we can take two different points $x_1$, $x_2$ in $\Om_1$ and we can construct the respective weight functions $\psi^1$ and $\psi^2$. For each $k=1,2$, $\psi^k$ is a transmission 
weight function in $\Om_0 \setminus \overline{B_\E(x_k)}$ and it remains to be shown that ($H_5$) is fulfilled in order to finish the proof of Proposition $\ref{psi12}$.

Let be $d = \frac{1}{2}|x_1 - x_2|$ such that $\E < d$. On the one hand, for all $x \in \Bi$ we have:
$$ \psi^1(x) \leq  \frac{a}{\rho_1^2} \E^2  + M \leq  \frac{a}{\alpha_1^2} \E^2  + M, $$
where $\alpha_1 = d(x_1,\Gm_1) > 0$.
On the other hand, if we denote $D_2 = \max\limits_{y \in \Gm_1}d(y,x_2)$, we get, for all $x \in \Bi$,
$$ \psi^2(x)  \geq \frac{a}{\rho_2^2} d^2  + M \geq  \frac{a}{D_2^2}d^2  + M.$$
Consequently, we have
 \begin{equation}
 \psi^2 -  \psi^1  \geq  a \left(  \frac{d^2}{D_2^2} -  \frac{\E^2 }{\alpha_1^2} \right) 
\hspace{0.5cm} \forall x \in B_{\E}(x_1).
\end{equation}
It is clear that an analogous result is true by interchanging $x_1$ and $x_2$ (now with $\alpha_2$ and $D_1$). 
Thus, taking $\E < \min\left( \frac{d \alpha_1}{D_2}, \frac{d \alpha_2}{D_1} \right) $ we 
get  ($H_5$)   and Proposition ~\ref{psi12} is proved.
\QED\\

>From Proposition \ref{psi12} and Theorem \ref{Carleman}, we obtain the following result:

\begin{Thm} \label{Carlemanconvex}
Let the coefficient $a$ be constant  in  the open set $\Om_j$   and equal to $a_j$ for
each $j = 1,2$.
Suppose that $~a_2 < a_1~$ and that  $\Om_1$ is an open,  bounded 
and strongly convex set with smooth boundary.
Then we have a Carleman estimate like $(\ref{ineqCarlem})$ for the Schr\"odinger  equation $(\ref{Anonlin})$ 
in the domain $\Om$.
\end{Thm}

\section{Stability of the inverse problem}\label{A3}

As described in the introduction,  will only give a local answer about the determination of the potential $p$. 
We will first work on a linearized version of the problem and consider the following Schr\"odinger equation :
\begin{equation}\label{Alin}
\left\lbrace
\begin{array}{ll}
iu'+~\Div\left(a(x)\nabla u\right)+q(x)u=f(x)R(x,t),~&\Omega\times(0,T)\\
u(x,t)=0,~&\partial\Omega \times(0,T)\\
u(x,0)=0,~&\Omega\\
\end{array} \right.
\end{equation}
Here we set $y=y(p)$ the weak solution to (\ref{Anonlin}) and $u=u(f)$
the one to (\ref{Alin}). If we formally linearize equation (\ref{Anonlin})
around a non stationary solution, we obtain equation (\ref{Alin}). In fact, we notice here
that if we set $f=p-q$, $u=y(q)-y(p)$ and $R=y(p)$ on $\Omega\times (0,T)$,
we obtain (\ref{Alin}) after substraction of (\ref{Anonlin}) with potential
$p$ from (\ref{Anonlin}) with potential $q$ and linearization.\\

\textbf{Linear inverse problem :}
Is it possible to determine $\left.f\right|_\Omega~$ from the knowledge of
the normal derivative $\left.\frac{\partial u}{\partial \nu} \right|_{\partial\Omega \times
(0,T)}$ where $R$ and $p$ are given and $u$ is the solution to (\ref{Alin})?\\

The following theorem proves that this inverse problem is well posed.

\begin{Thm}\label{ATh2}
Let $q\in L^{\infty}(\Omega)$ and $u$ be a solution of
equation $(\ref{Alin})$.
We assume that $$R \in W^{1,2}(0,T,L^{\infty}(\Omega)),$$
$$R(0)~~\textit{is real valued and}~~|R(x,0)|\geq r_0 > 0,~\textit{ a.e. in }~\overline{\Omega}.$$
There exists a constant $C=C(\Omega,T,\|q\|_{L^{\infty}(\Omega)},R)>0$
such that if
$$\dfrac{\partial u}{\partial\nu} \in H^1(0,T;L^2(\Gamma_0)),$$ then,
\begin{equation}\label{AT2}
\|f\|_{L^2(\Omega)}\leq C\left\|a_2\frac{\partial u}{\partial
\nu}\right\|_{H^1(0,T;L^2(\partial\Om))}.
\end{equation}
\end{Thm}
\noindent\textbf{Proof :}\\
As we need to estimate $\ds\frac{\partial u}{\partial \nu}$
in $H^1(0,T;L^2(\Gamma_0))$ norm, we work on the equation satisfied
by $v=u'$ :
\begin{equation}\label{Av}
\left\lbrace
\begin{array}{ll}
iv'+~\Div\left(a(x)\nabla v\right)+q(x)v=f(x)R'(x,t),~  &\Omega\times(0,T)\\
v(x,t)=0,~  &\partial\Omega \times(0,T)\\
v(x,0)=-if(x)R(x,0),~  &\Omega\\
\end{array} \right.
\end{equation}

The Carleman inequality we just obtained is the key of the proof. We extend the function
$v$ on $\Omega\times(-T,T)$ by the formula $v(x,t)=-\overline{v}(x,-t)$ for every 
$(x,t)\in\Omega\times(-T,0)$.  Since $R(0)$ and $f$ are real valued, $v \in C([-T,T];H_0^1(\Omega))$ 
and $\dfrac{\partial v}{\partial \nu}\in L^2((-T,T)\times\Gamma)$.
We also extend $~R~$ on $~\Omega\times(-T,T)$ by the formula $R(x,t)=\overline{R}(x,-t)$ for every $(x,t)\in\Omega\times(-T,0)$ 
and if we denote the extention of $R'$
by the same notation, then $R'\in L^2(-T,T;W^{1,\infty}(\Omega))$.
Thus, $v$ satisfies the same equation (\ref{Av}), set in $(-T,T)$.\\

As defined in Theorem \ref{Carleman},  for $k=1,2$, we set $w^k=e^{-s\varphi^k} v$ 
and $$P_1^{\psi^k}w^k=i\partial_tw^k+~\Div (a \nabla w^k)+s^2 a |\nabla\varphi^k|^2w^k.$$
Therefore, we define the following: 
$$I= \sum_{k=1}^{2}  Im \int_{-T}^0 \int _{\Omega} P_1^{\psi^k}w^k~\overline{w^k}\,dxdt.$$

On the one hand,
\begin{eqnarray}
I&=&\sum_{k=1}^{2}  Im \int_{-T}^0 \int _{\Omega} P_1^{\psi^k}w^k~\overline{w^k}\,dxdt\nonumber\\
&=&\sum_{k=1}^{2}  Im \int_{-T}^0 \int _{\Omega} \left(i\partial_tw^k+~\Div (a \nabla w^k)+s^2 a |\nabla\varphi^k|^2w^k\right)~\overline{w^k}\,dxdt \nonumber \\
&=&\sum_{k=1}^{2}  \int_{-T}^0 \int _{\Omega} Re\left(\partial_tw^k\overline{w^k}\right)
-Im\left(a \left|\nabla w^k\right|^2-s^2 a |\nabla\varphi^k|^2\left|w^k\right|^2\right)\,dxdt \nonumber\\
& = & \frac 12 \sum_{k=1}^{2} \int_{-T}^0 \int _{\Omega}\partial_t\left(|w^k|^2\right)\,dxdt \nonumber\\
&   =& \frac 12 \sum_{k=1}^{2} \int _{\Omega}|w^k(x,0)|^2\,dx \hspace{3cm}\nonumber\\
&   = &\frac 12 \sum_{k=1}^{2} \int _{\Omega}|f(x)|^2|R(x,0)|^2 e^{-2s\varphi^k(x,0)} \,dx.\label{exact}
\end{eqnarray}

On the other hand, Cauchy-Schwarz inequality and Carleman estimate
from  Theorem \ref{Carlemanconvex} give :
\begin{eqnarray*}
I&\leq& \sum_{k=1}^{2} \left( \int_{-T}^T \int_{\Omega}|P_1^{\psi^k}w^k|^2\,dxdt\right)^{\frac 12}
\left( \int_{-T}^T \int _{\Omega}|w^k|^2\,dxdt\right)^{\frac 12}\\
&\leq& \sum_{k=1}^{2}  \left\|P_1^{\psi^k}(w^k) \right \|_{L^2(Q)} \left\|w^k \right \|_{L^2(Q)}\\
&\leq& Cs^{-\frac32}  \sum_{k=1}^{2} \left( \left \|P^{\psi^k}(w^k) \right\|^2_{L^2(Q)}   
 +  s  \iint_{\Sigma_+^{\psi^k}} \theta^k \left|a \dn{w^k} \right|^2 \,d\sigma dt \right)\\
&\leq& Cs^{-\frac32}  \sum_{k=1}^{2} \left( \iint_{Q}|fR'|^2e^{-2s\varphi^k}\,dxdt 
+  s  \iint_{\Sigma_+^{\psi^k}}\theta^k \left|a_2 \dn{v} \right|^2e^{-2s\varphi^k} \,d\sigma dt \right).
\end{eqnarray*}
Then, $\varphi^k(x,t)=\frac{\alpha-e^{\lambda
\psi^k(x)}}{(T-t)(T+t)}$ is such that $e^{-2s\varphi^k(x,t)} \leq
e^{-2s\varphi^k(x,0)}$ for all $ x\in \Omega$ and $t\in (-T,T)$ and it
is easy to see that  $\theta e^{-2s\varphi}$ is bounded on $\Sigma_+^{\psi^k}$ 
and that using the definition of the extensions of $v$ and $R'$, we easily get
\begin{equation}
\label{majoration}
I  \leq  Cs^{-\frac32}  \sum_{k=1}^{2} \left( \int_{0}^{T}\hspace{-0,2cm}\int_{\Omega}|fR'|^2 e^{-2s\varphi^k(0)}\,dxdt 
+ s\iint_{\Sigma_+^{\psi^k}}\left|a_2\frac{\partial v}{\partial \nu}\right|^2\,d\sigma dt \right)
\end{equation}

>From $R \in W^{1,2}(0,T,L^{\infty}(\Omega))$ and $|R(x,0)|\geq r_0 > 0$ almost everywhere in $\overline{\Omega}$, we deduce that
$$\exists~g_0\in L^2(0,T),~ |R'(x,t)|\leq g_0(t)|R(x,0)|,~\forall
x\in \Omega,~t\in (0,T).$$
Hence, from (\ref{exact}) and (\ref{majoration}) we have :
\begin{eqnarray*}
&&\sum_{k=1}^{2} \int _{\Omega}|f|^2|R(0)|^2 e^{-2s\varphi^k(0)} \,dx \\
&\leq&Cs^{-\frac32}  \sum_{k=1}^{2} \left( \int_{0}^{T}\hspace{-0,2cm} \int_{\Omega}|fR'|^2 e^{-2s\varphi^k(0)}\,dxdt 
+ s\iint_{\Sigma_+^{\psi^k}}\left|a_2\frac{\partial v}{\partial \nu}\right|^2\,d\sigma dt \right)\\
&\leq& C  \sum_{k=1}^{2} s^{-\frac32} \int_{0}^{T}\hspace{-0,2cm} \int_{\Omega}|f|^2 |g_0|^2 |R(0)|^2 e^{-2s\varphi^k(0)}\,dxdt \\
&& +~C  \sum_{k=1}^{2} s^{-\frac12} \iint_{\Sigma_+^{\psi^k}}\left|a_2\frac{\partial v}{\partial \nu}\right|^2\,d\sigma dt.
\end{eqnarray*}
But $g_0\in L^2(0,T) $ implies $ \displaystyle\int_{0}^{T}|g_0(t)|^2\,dt \leq K <+\infty$
and so we can write
$$\left(1-\frac{CK}{s^{\frac32}}\right) \sum_{k=1}^{2} \int _{\Omega}|f|^2|R(0)|^2 e^{-2s\varphi^k(0)}\,dx
\leq C s^{-\frac12} \sum_{k=1}^{2}  \iint_{\Sigma_+^{\psi^k}}\left|a_2\frac{\partial v}{\partial \nu}\right|^2\,d\sigma dt$$
that becomes easily, if $s$ is large enough ($s>(CK)^{\frac 23}$) and $C$ remains a generic positive constant
$$\int _{\Omega}|f|^2|R(0)|^2 
\left( e^{-2s\varphi^1(0)} + e^{-2s\varphi^2(0)}\right)\,dx
\leq C s^{-\frac12} \iint_{\Sigma}\left|a_2\frac{\partial v}{\partial \nu}\right|^2\,d\sigma dt.$$
Moreover, since  $|R(x,0)|\geq r_0 > 0$  and 
$e^{-2s\varphi^k(x,0)} \geq e^{-2s \frac {\alpha-1}{T^2}}>0$ almost everywhere in $\overline{\Omega}$, we obtain
$$ \int _{\Omega}|f(x)|^2\,dx \leq C \iint_{\Sigma}\left|a_2\frac{\partial v}{\partial \nu}\right|^2\,d\sigma dt,$$\\
and therefore, Theorem \ref{ATh2} has been proved.
\QED\\

\noindent\textbf{Remark :} if we replace the assumption ~``$R(0) ~\textit{is real valued }$''~ by the following ~``$R(0)~\textit{takes its values in }~ i\mathbb R$''~,  then the appropriate extensions for $(x,t)$ in $\Omega\times(-T,0)$ are $v(x,t)=\overline{v}(x,-t)$ and $R(x,t)=-\overline{R}(x,-t)$. \\

 \smallskip

We will end this paper by the proof of Theorem~\ref{ip} which is a direct consequence of Theorem~\ref{ATh2}. Indeed, if  we set $\tilde u = y(q) - y(p)$, $f = p - q$ and $R = y(p)$, then $\tilde u$ is the solution of 
\begin{equation}\label{pbl}
\left\{ \begin{array}{ll}
\tilde u'  + \textnormal{div}( a \nabla \tilde u) +( p-f)\tilde u   =  f(x)R(x,t) & (0,T)\times\Omega\\
\tilde u  =  0 &(0,T)\times\Sigma \\
\tilde u(0) =0& \Om
\end{array}
\right. 
\end{equation}
where $q=p-f \in \mathcal U$, with $\mathcal U$ bounded in $L^\infty(\Om)$ from the hypothesis of Theorem~\ref{ip}. The key point is that in the proof of Theorem~\ref{ATh2}, all the constants $C>0$ depend on the $L^\infty$-norm of the potential. Thus, with $q\in \mathcal U$, we are actually, with equation (\ref{pbl}) in a situation similar to the linear inverse problem related to equation (\ref{Alin}) and we then obtain the desired result.
\QED

\section*{Appendix: Direct proof of Lemma \ref{positDef}}

Without lost of generality, we can take  $x_0 = 0$.
 Now, take $x , v \in \RR^N\setminus \{ 0\}$ and define  
 $g(t) = \mu^2(x+tv)$ for $t \in \RR$.
 Then $g$ depends only on the  restriction of $\mu^2$ to the plane  $ \Pi = \left< \{  x,v  \} \right> \subset \RR^N$ spanned by the vectors $x$ and $v$. Moreover, by definition of $\rho$, it is not difficult to see that $\rho|_\Pi = \rho_0$, where we have denoted by $\rho_0$  the function defined in the plane $\Pi$ as in $(\ref{defrho})$, 
but where    the closed  curve is given 
by $\Gamma_1 =  \Pi \cap \partial \Om_1$, wich by hypothesis is strongly convex (see Remark  \ref{planes}). 
 
It is not difficult to see  that  $\ds{ \frac{d^2g}{dt^2}(0)} = D^2(\mu^2)(x)(v,v)$ 
and then this expression depends only   on the curve $\Gamma_1 \subset \Pi$. We conclude that it suffices to consider   the two-dimensional case.

Assuming $N=2$,  $\Gamma_1$ can be parameterized in polar coordinates by 
 $$ \gamma(\theta) = (\rho(\theta)\cos\theta, \rho(\theta)\sin\theta) \quad  \theta \in [0, 2 \pi) .$$
The expression for the Hessian matrix of second derivatives in polar coordinates is
$$ D^2(\mu^2) = Q_\theta H(\mu^2) Q_\theta^T  $$
where $Q_\theta$ is the rotation matrix by angle $\theta$, and
\begin{equation}
H(\mu^2) = \left(  \begin{array}{cc}
\frac{\partial^2 \mu^2}{\partial r^2}    & \frac{1}{r}\left(\frac{\partial^2 \mu^2}{\partial r \partial \theta}  
     - \frac{1}{r} \frac{\partial \mu^2}{\partial \theta} \right) \vspace{0,3cm}\\
   \frac{1}{r}\left(\frac{\partial^2 \mu^2}{\partial r \partial \theta}  
     - \frac{1}{r} \frac{\partial \mu^2}{\partial \theta} \right)  & 
   \frac{1}{r^2} \frac{\partial^2 \mu^2}{\partial \theta^2} + \frac{1}{r} \frac{\partial \mu^2}{\partial r} 
\end{array}
\right). \nonumber
\end{equation}
Now, since $x_0 = 0$, we have $~ \mu^2(\theta, r) = \dfrac{\bar{a}}{\rho(\theta)^2} r^2  +M$.\\
One can notice that $\mu^2$ is well defined and smooth in $\Om_0 \setminus B_\E(x_0)$ (which means $\{r \geq \E\} \setminus  \Gamma_1$). All the computations that follows are valid in this set. We already said above that $\rho$ is constant with respect to $r$ and only depends on $\theta$ such that $\frac{\partial \rho}{\partial r} = 0$. Hence, we get 
\begin{equation}
\label{Hess}
H(\mu^2) = \frac{2 \bar{a}}{\rho^2}
\left(  \begin{array}{cc}
1   & -\frac{\rho_{\theta}}{\rho}  \\
-\frac{\rho_{\theta}}{\rho}   &  
\frac{1}{\rho^2}( 3\rho_{\theta}^2 - \rho\rho_{\theta \theta} + \rho^2)
\end{array}
\right),
\end{equation}
where we have denoted $\rho_{\theta} = \frac{\partial \rho}{\partial \theta}$ and so on.

We  will use the following well known lemma (see  \cite{Gray}) concerning curves in the plane:

\begin{Lem} Let $\gamma$ be a $C^2$ curve in the plane parameterized in polar coordinates by its angle: $\gamma(\theta) = (r(\theta)\cos\theta, r(\theta)\sin\theta)$. Then, the curvature of $\gamma$ is given by the formula
$$ \kappa_{\gamma}(\theta) = \frac{r^2 + 2r_{\theta}^2- r r_{\theta \theta}}{(r^2 + r_{\theta}^2)^{3/2}}.$$
\end{Lem}

Since the polar parametrization of $\Gm_1$ is given by $r(\theta) = \rho(\theta)$ and  $\Om_1$ is strongly convex, we obtain $$ \kappa_{_{\Gamma_1}} (\theta) = \frac{\rho^2 + 2\rho_{\theta}^2- \rho \rho_{\theta \theta}}{(\rho^2 + \rho_{\theta}^2)^{3/2}} > 0  \hspace{1cm} \forall \theta \in [0, 2 \pi[.$$
 
We will now prove that $H(\mu^2)$ is uniformly positive definite in $\Om_0 \setminus B_\E(x_0)$,  which will imply that the hypothesis $(H_4)$ is fulfilled in this set, since $a$ is piecewise constant. \\
 
 The eigenvalues of the matrix $\frac{\rho^2}{2 \bar{a}} H(\mu^2)$ satisfy the equation
$$r  =  \frac12 \left( d \pm \sqrt{d^2 - 4m} \right)$$
where $\ds{ d = \frac1{\rho^2} \left(    3\rho_{\theta}^2 - \rho\rho_{\theta \theta} + 2 \rho^2   \right) }$
 and
$$\ds{ m = \frac1{\rho^2} \left(  2 \rho_{\theta}^2 - \rho\rho_{\theta \theta} +  \rho^2   \right) 
   =      \frac{(\rho^2 + \rho_{\theta}^2)^{3/2}}{\rho^2}      \kappa_{_{\Gamma_1}} (\theta) > 0. }$$

Then $r_2 =  \dfrac12 \left( d + \sqrt{d^2 - 4m} \right) \le d$,
and $r_1 = \dfrac{m}{r_2} \ge \dfrac{m}{d} \quad \mbox{ for all } \theta \in [0, 2\pi).$
Since $\Om_1$ is bounded we get the desired result.  
\QED

\subsection*{Acknowledgments.}
The second author would like to acknowledge the partial support by  ECOS  grant  C04E08 and NSF grant DMS 0554571.


\medskip

\end{document}